\documentclass[10pt]{amsart}
\usepackage{tikz}
\usepackage{graphicx}
\usepackage{amssymb,mathrsfs}
\usepackage[colorlinks=true,
                    linkcolor=blue,
                    urlcolor=blue,
                    citecolor=blue,
                    anchorcolor=blue]{hyperref}
\usepackage{mathtools}
\usepackage{amsmath,amssymb}

\usepackage{rotating}
\usepackage[normalem]{ulem}

\theoremstyle{plain}
\newtheorem{theorem}{Theorem}[section]
\newtheorem{lemma}[theorem]{Lemma}
\newtheorem{proposition}[theorem]{Proposition}
\newtheorem{corollary}[theorem]{Corollary}

\theoremstyle{remark}

\theoremstyle{definition}

\newtheorem*{main}{Main Theorem}
\newtheorem*{mainproof}{Main Theorem Proof}

\def\Z{{\mathbb{Z}}}

\def\S{{\mathbb{S}}}
\def\Mod{{\rm Mod}}
\def\ScrC{{\mathscr{C}}}
\begin{document}

\title[Low-slope Lefschetz fibrations with finitely presented fundamental groups] 
{Low-slope Lefschetz fibrations with finitely presented fundamental groups}

\author[T{\"{u}}l\.{i}n Altun{\"{o}}z and Adalet \c{C}engel]
{T{\"{u}}l\.{i}n Altun{\"{o}}z and Adalet \c{C}engel}

\thanks{The second author was supported by Scientific and Technological Research Council of T\"urkiye (T\"UB\.{I}TAK) under Grant Number 124F502.}

\address{Faculty of Engineering, Ba\c{s}kent University, Ankara, Turkey}
\email{tulinaltunoz@baskent.edu.tr}
\address{Department of Mathematics, Bart\.{i}n University, Bart\.{i}n, Turkey}
\email{acengel@bartin.edu.tr}

\begin{abstract}
We prove that for any finitely presented group $G$, there exists a genus-$3k (k\in \Z)$ Lefschetz fibration over the two-sphere having slope as close as desired to $2$ and the fundamental group of the total space is isomorphic to given group $G$. These are the first examples not satisfying the slope inequality with total spaces being non-simply connected. As a corollary, by taking fiber sum of these Lefschetz fibrations with trivial genus-$3k$ Lefschetz fibration over any genus-$h(>0)$ (positive genus) surfaces, we obtain Lefschetz fibrations having slope arbitrarily close to $2$ with different fundamental groups over any genus-$h$ surface.
\end{abstract}

\maketitle

\section{Introduction} 

In algebraic geometry, the slope of a fibration is a  numerical invariant for surfaces fibered over curves. In their geography, it provides a significant relationship between the relative canonical class and the relative holomorphic Euler characteristic.

The notion of slope can be extended to four dimensional topology which is completely expressed in terms of the topological invariants of the total space. It also enables naturally to study the difference between smooth four-manifolds and algebraic (or complex) surfaces. More precisiously, it is a way to determine whether a Lefschetz fibration is holomorphic or not which is a complex geometric condition.

Let $f \colon X \to \S^2$ be a relatively minimal
genus-$g$ Lefschetz fibration where $g \geq 2$.
The slope of $f$ is defined by
\[
\lambda_f=\frac{K_f^2}{\chi_f},
\]
where
\[
K_f^2=c_1^2(X)+8(g-1)
\]
and
\[
\chi_f=\chi_h(X)+g-1.
\]
Here,
\[
c_1^2(X)=2e(X)+3\sigma(X)  \quad  \textrm{and} \quad 
\chi_h(X)=\frac{e(X)+\sigma(X)}{4},
\]
where $e(X)$ and $\sigma(X)$ denote the Euler
characteristic and signature of $X$, respectively.

Consequently, the slope can be expressed as
\[
\lambda_f=
\frac{
4\bigl(2e(X)+3\sigma(X)+8(g-1)\bigr)
}{
e(X)+\sigma(X)+4(g-1)
}.
\]

For a genus-$g$ Lefschetz fibration over $\S^2$
with $n$ critical points, we have
\[
e(X)=4-4g+n.
\]
Hence,
\[
\lambda_f=
4\frac{2n+3\sigma(X)}{n+\sigma(X)}.
\]

Xiao \cite{xiao} showed that the slope inequality $\lambda_f\geq 4-\frac{4}{g}$ is satisfied by every relatively minimal holomorphic non-trivial genus-$g$ Lefschetz fibrations over $\S^2$. This inequality gives a lower bound for the holomorphic counterparts of Lefschetz fibrations in the geography of algebraic fibrations. Unlike the holomorphic case, (symplectic) Lefschetz fibrations need not to satisfy the classical slope inequality \cite{monden,hkm,miyachishiga}. This distinction has led to
construction of Lefschetz fibrations with various
slope values and the investigation of how the slope
is influenced by the topology of the total space
and the monodromy of the fibration. Motivated by this, there is an analogous bound estimation for Lefschetz fibrations over $\S^2$ by Cengel and Korkmaz. In \cite{LowslopeLefschetz}, they construct examples of Lefschetz fibrations over $\S^2$ with slope converging to $2$ which is proposed as lower slope bound for Lefschetz fibrations over the two-sphere. Each fundamantel group of the corresponding total spaces are trivial.

Let $f\colon X\to\Sigma_h$ be a genus-$g$
Lefschetz fibration over a closed oriented
surface of genus $h$ where $g\geq 2$. The slope
of $f$ is defined as
\[
\lambda_f=\frac{K_f^2}{\chi_f}
\]
provided that $\chi_f\neq 0$, where
\[
K_f^2=2e(X)+3\sigma(X)-8(g-1)(h-1)
\]
and
\[
\chi_f=\frac{e(X)+\sigma(X)}{4}
       -(g-1)(h-1).
\]
Miyachi and Shiga \cite{miyachishiga} examine the necessary conditions for being holomorphic for Lefschetz fibrations over positive genus surfaces: the slope inequality and irreduciblity of the holonomoy group.  They show that neither of them imply the other one by giving counter examples. On the other hand, the holonomy group is irreducible for all Lefschetz fibrations over $\S^2$ so it does not give any interesting result for them. In \cite{miyachishiga}, Lefschetz fibrations over sufficiently large genus surfaces are given with negative slopes. This can not be extended to Lefschetz fibrations over $\S^2$. More concretely, the slope can not be negative because $K_f^2 \geq 4g-4$ and $\chi_f \geq 1.$ To the best of our knowledge the smallest slope is exactly $2$ which belongs to genus-$2$ hyperelliptic Lefschetz fibration.

A Lefschetz fibration over $\S^2$ can be described by its monodromy representation which is a positive factorization of identity in mapping class group. If it is over positive genus surface, there are some commutator factors are involved. There are two operations among many others to obtain new Lefschetz fibrations by using the existing ones: monodromy substitution and (twisted) fiber sum. The topological invariants of these constructed Lefschetz fibrations can be computed completely using 2-manifolds and related techniques. Even the Euler characteristic and signature may be uncontrolled during the operations, the slope of a fiber sum is a weighted average of the slopes of its summands. These constructions play an important role in the investigation of slope geography and its connections with the topology of symplectic four-manifolds. Moreover, examples with small slopes demonstrate significant differences between holomorphic and symplectic Lefschetz fibrations. These developments highlight the importance of understanding the range of possible slope values and their relationship with the topology and monodromy of Lefschetz fibrations.

Lefschetz fibrations in Main Theorem below are different from the rest by having total spaces with any finitely presented group. 

\begin{main}\label{main}
Let $G= < x_1, \ldots, x_k \mid R_1, \ldots, R_m >$ be a finitely presented group. For $k\geq 1$, there exists a sequence of genus-$3k$ Lefschetz fibrations over the two-sphere with slopes converge to $2$ and fundamental group of each one is isomorphic to $G$. 
\end{main}


 The group $G$ is a finitely presented group as in the main theorem. To prove the main theorem, we first construct genus-$3k$ Lefschetz fibrations over $\S^2$ with fundamental group is $G$ (Proposition~\ref{prop3.4}). We then perform the fiber sum and the subsequent chain substitution in such a way that the fundamental group of the resulting Lefschetz fibration is isomorphic to that of the
 initial Lefschetz fibration. More precisely, the gluing diffeomorphisms involved in the fiber sum are chosen such that the curves participating in the chain substitution are disjoint from those used to determine the generators and relations of the fundamental group. Consequently, chain substitution does not introduce additional relations affecting the fundamental group, and the resulting Lefschetz fibration has the same fundamental group as the initial one.
\qquad

\section{Preliminaries}
This section outlines the necessary background and gives the notation used throughout the paper.

Let $\Sigma_{g}^b$ be a compact, connected, oriented surface of genus $g$ with $b \geq 0$ boundary components. The mapping class group of this surface, denoted $\Mod(\Sigma_{g}^b)$, is defined as the group whose elements are isotopy classes of orientation-preserving self-diffeomorphisms of $\Sigma_{g}^b$, where every diffeomorphism is required to fix every point on the boundary. Also, the isotopies are also constrained to be the identity on the boundary. In the case of a surface with no boundary ($b=0$), we omit $b$ from the notation and write $\Sigma_{g}$ and $\Mod(\Sigma_{g})$.

Throughout this paper, a lowercase letter denotes a simple closed curve, and the corresponding capital letter denotes its positive (right-handed) Dehn twist. Our notation follows the functional order: for two diffeomorphisms $f$ and $g$, the composition $fg$ means $g$ first, then $f$. We denote the conjugation of $g$ by $f$ as $g^f$. All diffeomorphisms and curves are to be understood up to isotopy.

The mapping class group $\Mod(\Sigma_g^b)$ satisfies the following fundamental relations:
\begin{itemize}
    \item \textbf{Conjugation:} For any $f \in \Mod(\Sigma_{g}^b)$ and simple closed curves $x$ and $y$ with $f(x) = y$, the corresponding Dehn twists are conjugate:
        \[
        X^f = Y.
        \]
    \item \textbf{Commutativity:} If $x$ and $y$ are disjoint simple closed curves, then their Dehn twists commute:
        \[
        XY = YX.
        \]
    \item \textbf{Braid Relation:} If $x$ and $y$ intersect transversely at a single point, then their Dehn twists satisfy the braid relation:
        \[
        XYX= YXY.
        \]
         \begin{figure}[ht]
\begin{center}
\scalebox{0.3}{\includegraphics{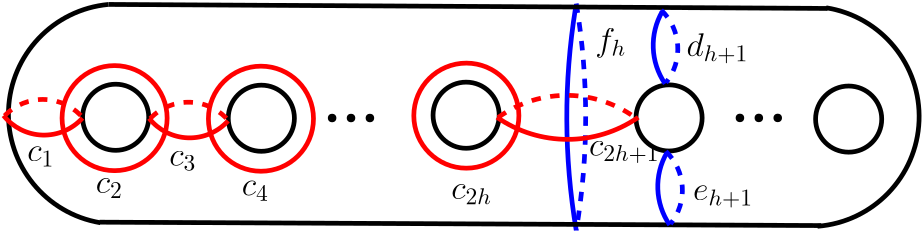}}
\caption{The curves in the chain relations.}
\label{chain}
\end{center}
\end{figure}
\item \textbf{Chain Relations:}
A chain of length $n$ on $\Sigma_g$ is an ordered
$n$-tuple $(c_1,\ldots,c_n)$ of simple closed curves
such that $c_i$ and $c_{i+1}$ intersect transversely
at exactly one point for $1\leq i\leq n-1$, while
$c_i$ and $c_j$ are disjoint whenever $|i-j|>1$.

If $n=2h+1$, a regular neighborhood of
$c_1\cup\cdots\cup c_{2h+1}$ is a genus-$h$
subsurface of $\Sigma_g$ with two boundary
components, denoted by $d_{h+1}$ and $e_{h+1}$,
as illustrated in Figure~\ref{chain}.
The corresponding \textit{odd chain relation} is
\[
(C_1 C_2 \cdots C_{2h+1})^{2h+2}
= D_{h+1}E_{h+1}.
\]
In this case, the element
\[
\ScrC_{2h+1}
=
(C_1 C_2 \cdots C_{2h+1})^{2h+2}
D_{h+1}^{-1}E_{h+1}^{-1}
\]
of $\Mod(\Sigma_g)$ is called the \textit{odd chain relator}.

If $n=2h$, a regular neighborhood of
$c_1\cup\cdots\cup c_{2h}$ is a genus-$h$
subsurface of $\Sigma_g$ with one boundary
component, denoted by $f_{h}$.
The \textit{even chain relation} is
\[
(C_1 C_2 \cdots C_{2h})^{4h+2}
=F_{h}.
\]
Then the corresponding relator is given by 
\[
\ScrC_{2h}
=
(C_1 C_2 \cdots C_{2h})^{4h+2}
F_{h}^{-1}
\]
of $\Mod(\Sigma_g)$ is called the \textit{even chain relator}.



\end{itemize}

\subsection{Lefschetz fibrations}

Let \( X \) be a closed oriented smooth 4-manifold and let $\mathbb{S}^2$ denote the two-sphere. A smooth map \( f: X \to \mathbb{S}^2 \) is called a \textit{genus-\(g\) Lefschetz fibration} on $X$ having finitely many critical values \( b_1, \ldots, b_n \) in such a way that there exists a unique critical point \( p_i \) in the singular fiber \( f^{-1}(b_i) \) such that around each \( p_i \) and \( b_i \) there are orientation--preserving local complex coordinate charts on which  \( f \) is given by \( f(z_1, z_2) = z_1^2z_2^2 \). (In general, the base space of a Lefschetz fibration can be any surface; however, we first restrict ourselves to the case where it is $\mathbb{S}^2$.) We assumed that \( f \) is relatively minimal. This means no fiber contains an exceptional sphere, i.e., a smoothly embedded sphere with self-intersection number $(-1)$. Each regular fiber of 
$f$ is a closed connected orientable genus-$g$ surface $\Sigma_g$; this positive integer $g$ is referring to as the genus of the Lefschetz fibration.

Each singular fiber of the Lefschetz fibration $(X,f)$ results from pinching  a specific simple closed curve, called a \textit{vanishing cycle} on a regular fiber to a single point. The vanishing cycle determines the diffeomorphism type of a regular neighborhood of the corresponding singular fiber.

For a genus-$g$ Lefschetz fibration $(X,f)$, let $b_0$ be a fixed regular value and let us identify the regular fiber $f^{-1}(b_0)$ with the genus-$g$ surface $\Sigma_g$. The Lefschetz fibration $(X,f)$ over $\mathbb{S}^2$ with vanishing cycles $\{ a_i\}$ is then encoded by a \textit{monodromy factorization}: a relation
\[
A_1 A_2 \cdots A_n = 1 \quad \text{in } \Mod(\Sigma_g).
\]

The monodromy factorization of a Lefschetz fibration is not unique, ; however, any two representations of which are related by two fundamental operations:

\textit{Hurwitz Moves:} This move swaps the order of Dehn twists $A_i$ and $A_{i+1}$ in the factorization, while conjugating one of them to preserve the overall product, which is defined as:
    \[
    A_1A_2 \cdots A_i A_{i+1} \cdots A_n \sim A_1A_2 \cdots \left(A_{i+1}\right)^{A_{i}}A_i\cdots A_n,
    \]
    where \(\left(A_{i+1}\right)^{A_{i}} = A_iA_{i+1} A_i^{-1}\) denotes the conjugation.
    
   \textit{Global Conjugation:} The entire factorization can be conjugated by a diffeomorphism \(\phi\) in $\Mod({\Sigma_g})$:
    \[
    A_1A_2 \cdots A_n \sim  C_1C_2 \cdots C_n,
    \]
    where \(\phi(a_i)=c_i\).

Let $f_i: X_i \to \mathbb{S}^2$ be a genus-$g$ Lefschetz fibration with monodromy factorization $W_i$, for each $i=1,2$. From each  $X_i$, we remove a  tubular neighborhood $\upsilon \Sigma_g$ of one chosen regular fiber 
$\Sigma_g$. Then we glue  the resulting manifolds along their boundaries using a specific diffeomorphism $\phi$, which yields a genus $g$ Lefschetz fibration
\[
f_1\sharp f_2:X_1 \sharp X_2 \to \mathbb{S}^2
\]
with monodromy factorization $W_1W_2^{\phi}$. This operation is called the \textit{(twisted) fiber sum} of $X_1$ and $X_2$.

For a genus-$g$ Lefschetz fibration $f:X\to \mathbb{S}^2$ with $n$ singular fibers, the Euler characteristic is given by $e(X)=4-4g+n$. Another topological invariant is the signature, $\sigma(X)$, which can be computed using various methods~\cite{Cengel2020,endo,endonagami,matsumoto2,ozbagci,smith}. 

Let $x_1, x_2, \dotsc, x_m$ and $y_1, y_2, \dotsc, y_n$ be simple closed curves on $\Sigma_g$ whose positive Dehn twists satisfy the relation
\[
X_1 X_2 \cdots X_m = Y_1 Y_2 \cdots Y_n.
\]

Given a positive relator $W = U \cdot X_1 X_2 \cdots X_m \cdot V$, we can form a new positive relator
\[
W' = U \cdot Y_1 Y_2 \cdots Y_n \cdot V
\]
by replacing the subword $X_1 X_2 \cdots X_m$ in $W$ with $Y_1 Y_2 \cdots Y_m$. In this case, letting $R = Y_1 Y_2 \cdots Y_n X_m^{-1} \cdots X_2^{-1} X_1^{-1}$, we say that $W'$ is obtained from $W$ by an \textit{$R$-substitution}. Thus, a new Lefschetz fibration with monodromy $W'$ can be constructed from a fibration with monodromy $W$ by performing an $R$-substitution.

Suppose that $1\leq h \leq g$. The odd chain relator $\ScrC_{2h+1}$ can be embedded in $\Mod(\Sigma_g)$ (in this embedding, the curves in the relation can be as in Figure~\ref{chain}). The following lemma describes the effect of a $\ScrC_{2h+1}$-substitution on the slope of a Lefschetz fibration.

\begin{lemma}\cite{Cengel2020}
Let $g \geq 2$, and let $f_1 \colon X_1 \to S^2$ be a
genus-$g$ Lefschetz fibration. Suppose that
$f_2 \colon X_2 \to S^2$ is obtained from $f_1$
by applying a $C_{2h+1}$-substitution. Then
the invariants satisfy
\[
K_{f_2}^2 = K_{f_1}^2 + 2h^2
\]
and
\[
\chi_{f_2} = \chi_{f_1} + \frac{h^2+h}{2}.
\]
\end{lemma}
Before proceeding, we state how the fiber sum construction affects the slope of the total spaces. This is given by the following result:
\begin{lemma}\cite{LowslopeLefschetz}\label{lemmafibersum}
Let $f_i: X \to S^2$ be Lefschetz fibration of genus $g \geq 2$ for each $i=1,2$. For any fiber sum $f = f_1 \#f_2$, the following  relations hold for its invariants:
\[
K^2_f = K^2_{f_1} + K^2_{f_2} \quad \text{and} \quad \chi_f = \chi_{f_1} + \chi_{f_2}.
\]
It follows immediately that the slope of the resulting fibration is
\[
\lambda_f = \frac{K^2_{f_1} + K^2_{f_2}}{\chi_{f_1} + \chi_{f_2}}.
\]
\end{lemma}
\section{Results}

For any positive integer $k$, Altun{\"{o}}z \cite{Altunoz2021} constructed genus-$3k$ Lefschetz fibrations $f_k^{(0)}: X_k^{(0)} \longrightarrow \S^2$, whose total spaces $X_k^{(0)}$ are diffeomorphic to $\Sigma_{k} \times \mathbb{S}^2 \#6\overline{\mathbb{CP}^2}$ and let $W_k$ denote the corresponding monodromy factorization. For the construction breeding of two generalized Matsumoto's genus-$2k$ Lefschetz fibrations are used. Now we use genus-$3k$ Lefschetz fibrations to obtain small slope Lefschetz fibrations with any finitely presented fundamental group. 
 \begin{figure}[ht]
\begin{center}
\scalebox{0.4}{\includegraphics{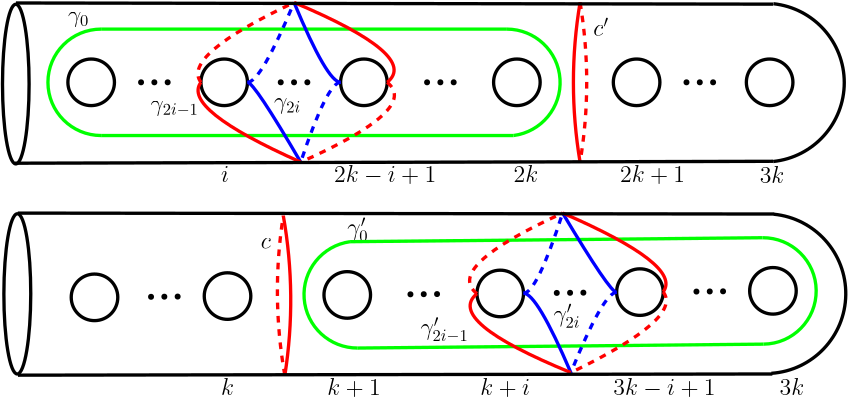}}
\caption{The curves in the monodromy of genus-$3k$ Lefschetz fibration.}
\label{genmatsumoto}
\end{center}
\end{figure}

\begin{figure}[ht]
\begin{center}
\scalebox{0.40}{\includegraphics{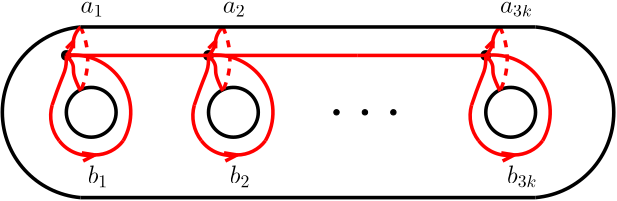}}
\caption{The generators of $\pi_1(\Sigma_{3k})$.}
\label{fundamental}
\end{center}
\end{figure}

In $\Mod(\Sigma_{3k}^1)$ 
\begin{equation}
\delta =
\Gamma_0\Gamma_1\Gamma_2\cdots \Gamma_{2k}C\Gamma_0 \Gamma_1\Gamma_2\cdots \Gamma_{2k}\Gamma_0^{'}\Gamma_1^{'}\Gamma_2^{'}\cdots \Gamma_{2k}^{'}C^{'}\Gamma_0^{'}\Gamma_1^{'}\Gamma_2^{'}\cdots \Gamma_{2k}^{'} \label{equation1}
\end{equation}
where the Dehn twist curves shown in Figure~\ref{genmatsumoto} gives the monodromy factorization of the Lefschetz fibration $X_k^{(0)} \to \mathbb{S}^2$. We want the boundary component of the surface to be on the left-hand side, which allows us to perform the desired substitution in the proof of main theorem. This can be achieved by using the breeding construction described in \cite{Altunoz2021}, in which the monodromy factorization with a boundary component is arranged so that the boundary lies on the left-hand side. Also, to avoid notational confusion, we rename the curves in  monodromy factorization given in~\cite{Altunoz2021}. We use capital letters to denote Dehn twists along the corresponding curves. Note that there are $8k+4$ nonseparating and two separating vanishing cycles in $W_k$. Since the fundamental group of $X_k^{(0)}$ is isomorphic to $\pi_1(\Sigma_{3k})$ divided out by the normal closure of the group generated by the vanishing cycles of $X_k^{(0)}$, 
\begin{equation}
\pi_1(X_k^{(0)})= \pi_1(\Sigma_{k})= < a_1, \ldots, a_k, b_1, \ldots, b_k \mid [a_1, b_1][a_2, b_2]\ldots [a_k, b_k]=1 >
\end{equation}
where the generators are as in Figure \ref{fundamental}.
\begin{theorem}\label{thm3.2}\cite{korkmaz2009}
Let $W=D_1D_2 \ldots D_n=1$ be a monodromy factorization of a Lefschetz fibration over the two-sphere with a section where its total space denoted by $X_W$. Let $e_1,e_2, \ldots, e_m$ be simple closed curves on fiber with the property that for each $i=1,2,\ldots,m,$ there exists an $j_i$ such that $e_i$ intersects $d_{j_i}$ transversely at only one point. By conjugating $W$ separately with positive Dehn twists along with each $e_m$'s followed by fiber sums give a new Lefschetz fibration. The fundamental group of the total space is
\[
\pi_1(X_W\#X_W^{E_1}\cdots\#X_W^{E_m})= \pi_1(X_W) \big/<e_1,e_2, \ldots, e_m>.
 \]
\end{theorem}
For the proof of this theorem, see Korkmaz \cite[Theorem $3.2$]{korkmaz2009}.

\begin{corollary}
For each $k\geq 1$, there exists a genus-$3k$ Lefschetz fibration $f_k^{(1)}:X_k^{(1)} \to \mathbb{S}^2$
with fundamental group isomorphic to a free group on $k$ generators.
\end{corollary}

\begin{proof}
By applying Theorem \ref{thm3.2} to a Lefschetz fibration with monodromy as in  (\ref{equation1}), we obtain a new Lefschetz fibration which has monodromy factorization
\begin{equation}
  W_k W_{k}^{A_1} W_{k}^{A_2}\cdots W_{k}^{A_k}  \label{equation3}.
\end{equation}

Let us denote the total space of the genus-$3k$ Lefschetz fibration with monodromy (\ref{equation3}) by $X_k^{(1)}$. Then the fundamental group $\pi_1(X_k^{(1)})$ is a free group isomorphic to $F_k$ and represented by   

\begin{equation}
\pi_1(X_k^{(1)})= < b_1, \ldots, b_k  >.\label{equation4}
\end{equation}
\end{proof}
Consider a finitely presented group $G = < x_1, \ldots, x_k \mid R_1, \ldots, R_m > $. It is shown that for a given finitely presented group $G$, it is possible to construct a Lefschetz fibration with fundamental group isomorphic to $G$ \cite{korkmaz2009, kobayashimonden}. We will use the following lemma to construct Lefschetz fibrations with fundamental group given G.

\begin{lemma}\cite{GGP}\label{3.3}
Let $G$ be a finitely presented group. Then it has a presentation $ < x_1, \ldots, x_k \mid R_1, \ldots, R_m >$ with the following properties: 

\begin{enumerate}
    \item each $R_j$ is a positive (no inverses) word in $ x_1, \ldots, x_k $;
    \item each generator $x_i$ appears at most once in each $R_j$;
     \item the cyclic order (by index) of the generators $ x_1, \ldots, x_k $ is preserved in each $R_j$.
\end{enumerate}
\end{lemma}

This lemma is important because this kind of group representations make possible the correspondence between the relations in the group and vanishing cycles of Lefschetz fibration without any self intersection. In the proof of the following proposition, we follow the approach used in the proof of Theorem $A$ in~\cite{arabadjibaykur}.

\begin{proposition} \label{prop3.4}
Let $G$ be a finitely presented group with properties as in Lemma \ref{3.3}. For each $k\geq 1$, there exists a genus-$3k$ Lefschetz fibration $f_k^{(2)}:X_k^{(2)} \to \mathbb{S}^2$
with fundamental group isomorphic to $G$.
\end{proposition}

\begin{figure}[ht]
\begin{center}
\scalebox{0.25}{\includegraphics{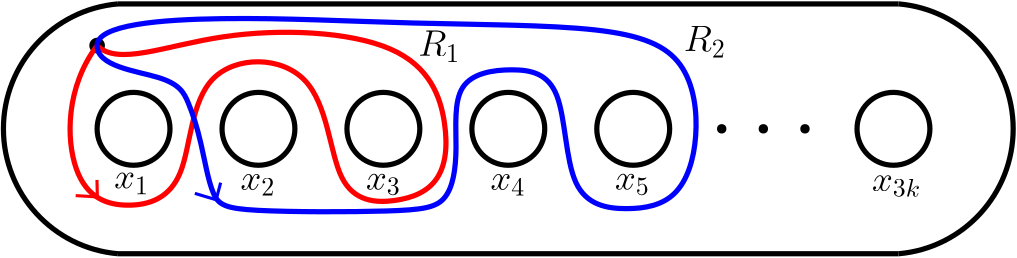}}
\caption{Some examples of relator curves on $\Sigma_{3k}$.}
\label{relators}
\end{center}
\end{figure}

\begin{proof}
Consider the Lefschetz fibration with monodromy as in (\ref{equation3}) such that the curves $b_i$'s are the generators of fundamental group of the corresponding
total space. Then identify each $b_i$ with given group generators $x_i$. Note that the number of $k$ generators in $G$ determines the genus of the regular fiber surface, $3k$. Assume that all the relators in $G$ can be nicely represented by the embedded curves as in Figure~\ref{relators} in terms of $x_i$'s where for example $R_1=x_1x_3, R_2=x_2x_3x_5$ and so on.
As a result, we construct a Lefschetz fibration $f_k^{(2)}:X_k^{(2)}\to \mathbb{S}^2$ with monodromy
\begin{equation}
 V= W_k W_{k}^{A_1} W_{k}^{A_2}\cdots W_{k}^{A_k} W_{k}^{R_1}W_{k}^{R_2}\cdots W_{k}^{R_m}=1. \label{equation5}
\end{equation}
The fundamental group of $X_k^{(2)}$ is isomorphic to
\begin{equation}
\pi_1(X_k^{(2)})= < x_1, \ldots, x_k \mid R_1, R_2,\ldots,R_m > =G.\label{equation6}
\end{equation}
\end{proof}
Now we give a result on small slopes that we use in our construction. The proof and the details can be found in \cite[Theorem $4.1$]{LowslopeLefschetz}.

\begin{theorem}\cite{LowslopeLefschetz}\label{thm4.1}
Let $g\geq3, r\geq1$ and $1 \leq h \leq g-1$ be integers and let $f_0:X_{0}\to \mathbb{S}^2$ be a genus-$g$ Lefschetz fibration with the monodromy factorization $VC^r$, where $C$ is a positive Dehn twist about a nonseparating curve and $V$ is a product of positive Dehn twists. For each integer $n\geq1$, there is a genus-$g$ Lefschetz fibration $f_n:X_{n}\to \mathbb{S}^2$ such that

\begin{enumerate}
    \item $K_{f_n}^2= 2^nK_{f_0}^2+2^nr[(h+1)^n-1]h$
    \item $\chi_{f_n}= 2^n\chi_{f_0}+2^{n-2}r[(h+1)^n-1](h+1)$.
\end{enumerate}
In particular,the slope of $f_n$ is 
\[
\lambda_{f_n}=\frac{4K_{f_0}^2\,+\, 4r[(h+1)^n-1]h}{4\chi_{f_0}\,+\,r[(h+1)^n-1](h+1)}
\]
so that its limit is $\frac{4h}{h+1}$.
\end{theorem}

In this theorem, the value of slope in result is independent of the initial choice of Lefschetz fibration. It is proved that the theorem starting with hyperelliptic Lefschetz fibration whose slope is exactly on the bound $4-4/g$. Then they claimed that whether $\lambda_{f_0}>\frac{4h}{h+1}$ or $\lambda_{f_0}<\frac{4h}{h+1}$, by applying the precedure iteratively, $\lambda_{f_n}$ converges to $\frac{4h}{h+1}$. To the best of our knowledge we don't know any Lefschetz fibration with slope less than $\frac{4h}{h+1}$ when $h=1$.
\begin{lemma}\label{lemma $3.6$}
There exists a sequence of genus-$3k$ Lefschetz fibrations over $\mathbb{S}^2$ with slope converging to $2$.
\end{lemma}
\begin{proof}


We follow the approach used in the proof of Theorem~\ref{thm4.1}, applying suitable chain substitutions for \(h=1\).
Since the topological type of the curves $\gamma_{2k}^{'}, d_{3k-1}$ and $e_{3k-1}$ are the same, there exist diffeomorphisms $\phi_i^1(i=1,2)$ sending  $\gamma_{2k}^{'}$ to $d_{3k-1}$ and $e_{3k-1}$ respectively. One can consider the following ones for $\phi_i^1(i=1,2)$:
\[
\phi_1^1=C_{6k-2}D_{3k-1}T \quad  \text{and} \quad \phi_2^1=C_{6k-2}E_{3k-1}T,
\]
where $T=C_{6k-3}C_{6k-2}C_{6k-4}C_{6k-3}\cdots C_{4k+1}C_{4k+2}\Gamma_{2k}^{'}C_{4k+1}$ (the Dehn twist curves in $\phi_i^1(i=1,2)$ are depicted in Figure \ref{phi}).

\begin{figure}[ht]
\begin{center}
\scalebox{0.3}{\includegraphics{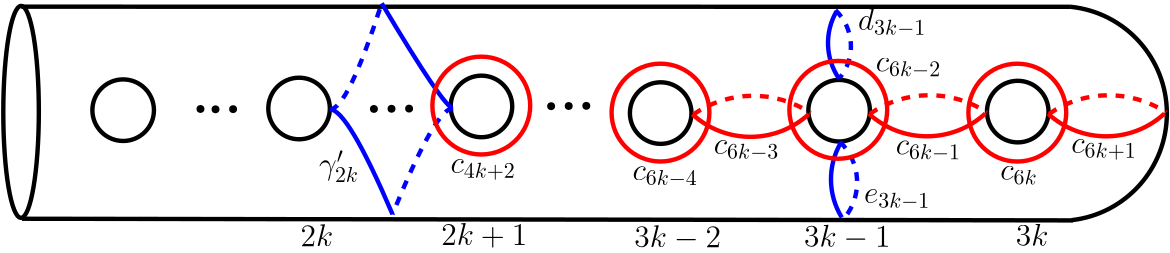}}
\caption{The Dehn twist curves curves in diffeomorphisms $\phi_i^1(i=1,2)$  on $\Sigma_{3k}$.}
\label{phi}
\end{center}
\end{figure}
Then conjugations of the monodromy in (\ref{equation5}) with $\phi_i^1(i=1,2)$ followed by a fiber sum and Hurwitz moves give
\begin{align}
  V_{0,0} &= V^{\phi_1^1}V^{\phi_2^1} \nonumber \\
    &= P_0^{\phi_1^1}(\Gamma_2^{'})^{\phi_1^1}P_0^{\phi_2^1}(\Gamma_2^{'})^{\phi_2^1} \nonumber \\
    &=\Phi_0 D_{3k-1}E_{3k-1} \label{eq6}
\end{align}
where $\Phi_0$ is a product of $(16k+10)$ positive Dehn twists. Then applying a chain substitution of lenght $3$
\begin{align}
\Delta_3:=D_{3k-1}E_{3k-1}=(C_{6k+1}C_{6k}C_{6k-1})^4 \label{eq7}
\end{align}
to monodromy in (\ref{eq6}), we have a new factorization 
\begin{align}
  V_{0,1} &= \Phi_0 (C_{6k+1}C_{6k}C_{6k-1})^4  \label{eq8}
\end{align}
which corresponds to a new Lefschetz fibration $f_{k,1}:X_{k,1}\to \mathbb{S}^2$ with slope $\lambda_{f_{k,1}} = 8-\frac{18}{(4k+1)}.$  By applying the same idea, that is choosing a diffeomorphism 
$\phi_i^2(i=1,2)$ sending the curve $c_{6k+1}$ to $d_{3k-1}$ and $e_{3k-1}$ respectively. Here one can choose these difeeomorphisms as follows.
\[
\phi_1^2=C_{6k-2}D_{3k-1}C_{6k-1}C_{6k-2} \quad  \text{and} \quad \phi_2^2=C_{6k-2}E_{3k-1}C_{6k-1}C_{6k-2},
\]
where the Dehn twist curves as shown in Figure~\ref{phi}.
Then we may write
\begin{align}
  V_{1,0} &= V_{0,1}^{\phi_1^2}V_{0,1}^{\phi_2^2} \nonumber \\
    &= P_1^{\phi_1^2}(C_{6k+1})^{\phi_1^2}P_1^{\phi_2^2}(C_{6k+1})^{\phi_2^2} \nonumber \\
    &=\Phi_1 D_{3k-1}E_{3k-1} \label{eq9}
\end{align}
where $\Phi_1$ is a product of $(32k+42)$ positive Dehn twists. Substituting $\Delta_3$ gives
\begin{align}
  V_{1,1} &= \Phi_1 (C_{6k+1}C_{6k}C_{6k-1})^4  \label{eq10}
\end{align}
where the corresponding Lefschetz fibration $f_{k,2}:X_{k,2}\to \mathbb{S}^2$. 
The monodromy (\ref{eq10}) again contains the curve $c_{6k+1}$, using the diffeomorphisms $\phi_i^3(i=1,2)$ sending the curve $c_{6k+1}$ to $d_{3k-1}$ and $e_{3k-1}$ respectively, we get a new monodromy $V_{2,0}$.
One can iterate this as arbitrarily as many times. Therefore, for each integer $i\geq1$, we get the monodromy
\begin{align}
  V_{i,0} &= V_{{i-1},{1}}^{\phi_1^{i+1}}V_{{i-1},{1}}^{\phi_2^{i+1}} \nonumber \\
    &= P_i^{\phi_1^{i+1}}(C_{6k+1})^{\phi_1^{i+1}}P_i^{\phi_2^{i+1}}(C_{6k+1})^{\phi_2^{i+1}} \nonumber \\
    &=\Phi_i D_{3k-1}E_{3k-1} \label{eq11}
\end{align}
and then the monodromy
\begin{align}
  V_{i,1} &= \Phi_i (C_{6k+1}C_{6k}C_{6k-1})^4. \label{eq12}
\end{align}
And there exists a genus-$3k$ Lefschetz fibration $f_{k,i}:X_{k,i}\to \mathbb{S}^2$  with monodromy factorization (\ref{eq10}). By Theorem~\ref{thm4.1}, the limit of the slope $\lambda_{k,i}$ of $f_{k,i}$ is equal to $2$. This finishes the proof.
\end{proof}

Now we are ready to give the proof the main theorem.

\begin{mainproof}\label{mainproof} Let $G$ be a finitely presented group having the presentation $ < x_1, \ldots, x_k \mid R_1, \ldots, R_m >$ as in Lemma~\ref{3.3}. For $k\geq 1$, consider the family of genus-$3k$
Lefschetz fibrations
$f_{k,i}:X_{k,i}\to \mathbb{S}^2$,
with monodromy factorizations given by \eqref{eq12}. The slopes of these fibrations converge to $2$
as $i\to\infty$. To construct these fibrations, we start with the genus-$3k$ Lefschetz fibration
$f_k^{(2)}: X_k^{(2)} \longrightarrow \S^2$ with monodromy factorization
\eqref{equation5}, whose fundamental group
is isomorphic to $G$ by Proposition~\ref{prop3.4}. Let $\alpha_1,\ldots,\alpha_n$ be the vanishing
cycles of the Lefschetz fibration $f_k^{(2)}: X_k^{(2)} \longrightarrow \S^2$, and let
$\beta_1,\ldots,\beta_r$ denote the vanishing
cycles of the fibration used in the fiber sum,
including those arising from chain substitutions. For $i,j=1,2$, the diffeomorphisms $\phi_i^j$ are the gluing diffeomorphism used to perform the fiber sum before applying the chain substitutions, then
\[
\pi_1 (X_{k,i})\cong
\frac{\pi_1(\Sigma_{3k})}
{\langle\
\alpha_1,\ldots,\alpha_n,
\phi_{i_1}^{j_1}(\beta_1),\ldots,
\phi_{i_r}^{j_r}(\beta_r)
\rangle},
\]
where $i_\ell,j_\ell\in\{1,2\}$ for
$\ell=1,\ldots,r$.
Since each $\phi_{i_\ell}^{j_\ell}(\beta_\ell)$,
for $\ell=1,\ldots,r$, can be expressed in terms
of generators of $\pi_1(\Sigma_{3k})$ that are
distinct from $x_1,\ldots,x_k$ and, all these elements represent the
identity in $\pi_1(X_{k,i})$. Therefore,  the additional relations are
redundant. Consequently,
\[
\pi_1(X_{k,i})\cong\pi_1( X_k^{(2)} )\cong G.
\]
This completes the proof.

\end{mainproof}

\subsection{Lefschetz fibrations over any genus-$h>0$ surface.}
In this section, we obtain Lefschetz fibrations over a closed orientable surface of $\Sigma_h$ of genus $h>0$ whose slope arbitrarily close to $2$ with non-simply connected total spaces. Our result is as follows.
\begin{corollary}
Let $G= < x_1, \ldots, x_k \mid R_1, \ldots, R_m >$ be a finitely presented group. For $k\geq 1$, there exists a sequence of genus-$3k$ Lefschetz fibrations over the surface $\Sigma_h$ of genus $h>0$  with slopes converge to $2$ and fundamental group of each one is isomorphic to the group $G\times \pi_1(\Sigma_h)$.     
\end{corollary}
\begin{proof}
 Let $f_{k,i}^{(h)}:X_{k,i}^{(h)}\longrightarrow\Sigma_h$ be the genus-$3k$ Lefschetz fibration over the surface $\Sigma_h$ obtained by taking fiber sum of genus-$3k$ Lefschetz fibration $f_{k,i}:X_{k,i}\longrightarrow \S^2$ constructed in the proof of main theorem and trivial genus-$3k$ Lefschetz fibration over the surface $\Sigma_h$. We first show that for each $k\geq 1$, the slope of $f_{k,i}^{(h)}:X_{k,i}^{(h)}\longrightarrow\Sigma_h$ converges to $2$. For $g=3k$,
 \[ e(X_{k,i}^{(h)})=e(X_{k,i})+4h(g-1), \qquad \sigma(X_{k,i}^{(h)})=\sigma(X_{k,i}), \]
 which gives
\[
\begin{aligned}
K_{f_{k,i}^{(h)}}^2
&=3\sigma(X_{k,i}+2[e(X_{k,i}+4h(g-1)]-8(g-1)(h-1)\\
&=3\sigma(X_{k,i})+2e(X_{k,i})+8(g-1)=K_{f_{k,i}}^2.
\end{aligned}
\]
Similarly,
\[
\begin{aligned}
\chi_{f_{k,i}^{(h)}}
&=\frac{\sigma(X_{k,i})+e(X_{k,i})+4h(g-1)}4-(g-1)(h-1)\\
&=\frac{\sigma(X_{k,i})+e(X_{k,i})}{4}+g-1=\chi_{f_{k,i}}.
\end{aligned}
\] Hence, $\lambda_{f_{k,i}^{(h)}}=\lambda_{f_{k,i}}$.

 For each positive integers $i$ and $k$, since the genus-$3k$ Lefschetz fibration $f_{k,i}:X_{k,i}\longrightarrow \S^2$ has a section, then the fiber sum along regular fibers yields
\[
\pi_1(X_{k,i}^{(h)})\cong G\times\pi_1(\Sigma_h), 
\] which completes the proof.
\end{proof}
\qquad
\qquad


\textbf{Acknowledgements} The authors would like to thank Mustafa Korkmaz for valuable feedback and encouraging comments on the results of this paper.
The second author of this paper was supported by Scientific and Technological Research Council of T\"urkiye (T\"UB\.{I}TAK) under Grant Number 124F502. The authors thank T\"UB\.{I}TAK for their support.

\section*{Declaration of AI Assistance}

The authors used ChatGPT (OpenAI) solely for
grammatical and linguistic improvements and
assistance with LaTeX formatting.
All mathematical results, proofs, and scientific
contributions presented in this paper are
entirely the authors' own.

\end{document}